\documentclass[a4paper]{article}
\usepackage{latexsym}
\usepackage{amssymb}

\newtheorem{Theorem} {Theorem} [section]
\newtheorem{Proposition} [Theorem] {Proposition}
\newcommand{\Proof}{ \noindent{\bf Proof.}\quad }
\newcommand{\qed}{\hfill$\Box$\medskip}

\newcommand{\adj}{\sim}

\newcommand{\scaps}[1]{\mbox{\rm\textsc{#1}}}

\title{The equivalence of two inequalities for quasisymmetric designs}
\author{A. E. Brouwer}
\date{2020-07-25}

\begin{document}
\maketitle

% \section{Introduction}
It has been an open problem whether Hobart's inequality on the
parameters of a quasisymmetric 2-design is independent of earlier
known restrictions. In this note we show that it is equivalent
to inequalities found by Neumaier and Calderbank.
We also give some more parameter sets ruled out by the Blokhuis-Calderbank
inequality.

\section{Quasisymmetric designs}
A {\em design} is a finite set called the point set,
provided with a collection of subsets called blocks.
A $t$-$(v,k,\lambda)$ {\em design} is a design with $v$ points,
where all blocks have size $k$ and any $t$ distinct points are
in precisely $\lambda$ blocks.

A {\em quasisymmetric design} with intersection numbers $x,y$,
is a design where distinct blocks meet in either $x$ or $y$ points,
where $x,y$ are distinct and both occur.

\medskip
A {\em strongly regular graph} with parameters $(v,k,\lambda,\mu)$
is a finite undirected graph without loops, having both edges and nonedges,
with $v$ vertices, regular of valency $k$, where two distinct adjacent
(resp. nonadjacent) vertices have precisely $\lambda$ (resp. $\mu$)
common neighbours. In this note we shall write $(V,K,\Lambda,M)$
for the parameters of a strongly regular graph, to avoid a clash with
design parameters.

\medskip
Let $(X,{\cal B})$ be a quasisymmetric 2-$(v,k,\lambda)$ design
with intersection numbers $x,y$, where $1 < k < v$.
The number of blocks on each point is $r = \lambda(v-1)/(k-1)$
and the total number of blocks is $b = vr/k$.

\medskip
Let $N$ be the point-block incidence matrix. Let $A$ be the 0-1 matrix
indexed by the blocks with $(B,C)$-entry 1 precisely when $|B \cap C| = x$.
Then $NN^\top = rI+\lambda (J-I)$ and $N^\top N = kI+xA+y(J-I-A)$.
Now $A$ is the adjacency matrix of a strongly regular graph.
Indeed, $NN^\top$ has two eigenvalues $r-\lambda$ and $kr$, so
$N^\top N$ has three eigenvalues $0$, $r-\lambda$ and $kr$, and also
$A = \frac{1}{x-y}(N^\top N - (k-y)I - yJ)$ has three eigenvalues,
namely $K = \frac{(r-1)k-(b-1)y}{x-y}$, $R = \frac{r-\lambda-k+y}{x-y}$
and $S = -\frac{k-y}{x-y}$ with multiplicities $1$, $v-1$, and $b-v$,
respectively.

\medskip
We see that the intersction-$x$ graph of $(X,{\cal B})$
with vertex set ${\cal B}$, where $B \adj C$ when $|B \cap C| = x$,
is strongly regular with parameters $(V,K,\Lambda,M)$ and eigenvalues
$K$, $R$, $S$, where $V = b$, and $K,R,S$ are as above, and $\Lambda,M$ are
determined by $RS = M-K$ and $R+S = \Lambda-M$.

% \medskip
% Quasi-symmetric designs were introduced by
% \scaps{Goethals} \& \scaps{Seidel} \cite{GoethalsSeidel69,GoethalsSeidel70}.
%
%% It is rumored that also Shrikhande & Bhagwandas,
%% Duals of incomplete block designs, J. Indian Statist. Assoc. 3 (1965) 30–37.
%% proved that a qs design gives a srg; I have not found this paper yet.
%

\medskip
Many examples are known. For example, the Steiner system
$S(4,7,23)$ is a quasisymmetric $2$-$(23,7,21)$ design
with intersection numbers 1 and 3. Its intersection-3 graph
is strongly regular with parameters $(V,K,\Lambda,M) = (253,140,87,65)$
with spectrum $140^1 ~ 25^{22} ~ (-3)^{230}$ where multiplicities
are written as exponents.

\medskip
\scaps{Blokhuis} \& \scaps{Haemers} \cite{BlokhuisHaemers01}
constructed an infinite family of examples with parameters
$v=q^3$, $k=\frac12 q^2 (q-1)$, $\lambda = \frac14 q(q^3-q^2-2)$,
$x = \frac12 k$, $y = x-\frac14 q^2$ where $q$ is a power of two.

\subsection{Complement}
Given a quasisymmetric 2-$(v,k,\lambda)$ design $(X,{\cal B})$,
with $b$ blocks, $r$ on each point, and intersection numbers
$x,y$, the {\em complementary design} is $(X,{\cal B}')$, where
${\cal B}' = \{X \setminus B \mid B \in {\cal B}\}$.
It has parameters $v' = v$, $k' = v-k$, $\lambda' = b-2r+\lambda$,
$b' = b$, $r' = b-r$, $x' = v-2k+x$, $y' = v-2k+y$.

\section{Inequalities}
% We shall need several known inequalities for the parameters
% of a quasisymmetric 2-design.

\subsection{The Calderbank-Cowen inequality}
The following result allows one to express the number of blocks $b$
of a quasi-symmetric 2-design in terms of the parameters $v,k,x,y$.

\begin{Proposition}\label{CCowen-ineq}
{\rm (\scaps{Calderbank} \cite{Calderbank88a})}
Every $1$-$(v,k,r)$ design with $b$ blocks, and two block intersection
numbers $x,y$, satisfies
$$
1 - \frac{1}{b} \le \frac{k(v-k)}{v(v-1)}
\left( \frac{(v-1)(2k-x-y)-k(v-k)}{(k-x)(k-y)}\right)
$$
with equality if and only if the design is a $2$-design. \qed
\end{Proposition}

\subsection{Neumaier's inequality}
Let $\Gamma$ be a strongly regular graph.
A proper nonempty subset $Y$ of its vertex set is called a {\em regular set}
with {\em degree $d$} and {\em nexus $e$} when each vertex inside
(resp. outside) $Y$ has $d$ (resp. $e$) neighbours in $Y$.

Let $\Gamma$ be the strongly regular graph on the blocks of a quasi-symmetric
$2$-$(v,k,\lambda)$ design $(X,{\cal B})$ with block intersection numbers
$x,y$, where blocks are adjacent if they meet in $x$ points.
Let $r = \lambda(v-1)/(k-1)$ be the replication number
(number of blocks on any point).

\begin{Proposition} \label{nexus}
{\rm (\scaps{Neumaier} \cite{Neumaier82a})}
The sets of all blocks $S(u)$ containing a fixed point $u$
are regular sets in $\Gamma$ of size $r$,
degree $d = \frac{(\lambda-1)(k-1)-(r-1)(y-1)}{x-y}$
and nexus $e = \frac{\lambda k - ry}{x-y}$.
\end{Proposition}
\Proof
Clearly, $|S(u)| = r$.
For $B \in S(u)$, with $d_B$ neighbours in $S(u)$,
count the number of pairs $(v,C)$ with $v \ne u$ and $C \ne B$ and $u,v \in C$
and $v \in B$.
This number is $(k-1)(\lambda-1)$ and also $d_B(x-1) + (r-d_B-1)(y-1)$
so that $d = d_B$ does not depend on $B$ and has the stated value.
Similarly, for $B \not\in S(u)$, with $e_B$ neighbours in $S(u)$,
we find $k \lambda = e_B x + (r-e_B)y$, so that $e_B$ does not depend on $B$
and has the stated value.
\qed

\begin{Proposition} \label{neum-ineq-qs}
{\rm (\scaps{Neumaier} \cite{Neumaier82a})}
The parameters of $(X,{\cal B})$ satisfy
$$
B(B-A) \le AC,\eqno{(N)}
$$
where
$$A = (v-1)(v-2),~~B = r(k-1)(k-2)$$
$$C = rd(x-1)(x-2) + r(r-1-d)(y-1)(y-2).$$
Equality holds if and only if $(X,{\cal B})$ is a 3-design.
\end{Proposition}
\Proof
For distinct points $u,v,w$, let $\lambda_{uvw}$ denote the number
of blocks containing these three points. Fix $u$ and sum over all
ordered pairs $v,w$ with $u,v,w$ distinct. One obtains
$\sum 1 = A$, $\sum \lambda_{uvw} = B$,
$\sum \lambda_{uvw}(\lambda_{uvw} - 1) = C$.
Now $0 \le \sum (\lambda_{uvw} - \frac{B}{A})^2 = B+C-\frac{B^2}{A}$.
\qed

%
% Sometimes rounding helps here
% v=839 k=343 lb=58653 x=147 y=133 has equality, but A does not divide B
% Here b=v(v-1)/2: r=143717, b=351541, N=0
% ./qsi 839 343 58653 147 133.
%
% v=1717 k=442 lb=5586 x=119 y=102 has equality, but A does not divide B
% ./qsi 1717 442 5586 119 102
%
% So far we do not have examples where rounding helps and a 3-design
% was possible.
%
%
% \medskip
% In rare cases one can strengthen this inequality by rounding:
% $$
% 0 \le \sum (\lambda_{uvw} - \lfloor\frac{B}{A}\rfloor)(\lambda_{uvw} - \lceil\frac{B}{A}\rceil)
% = B + C -B(\lfloor\frac{B}{A}\rfloor + \lceil\frac{B}{A}\rceil)
% + A \lfloor\frac{B}{A}\rfloor \lceil\frac{B}{A}\rceil .
% $$
% If $\frac{B}{A} = \lfloor\frac{B}{A}\rfloor + \eps$, with $0 < \eps < 1$, then
% the right-hand side of this latter inequality is $B+C-\frac{B^2}{A}-\eps(1-\eps)A$.
% This rules out e.g. $(v,k,\lambda,x,y) = (839,343,58653,147,133),
% (1717,442,5586,119,102)$ that have equality in (N) (but cannot belong
% to 3-designs).
%

One may check that Neumaier's inequality (N) for a design is equivalent
to the inequality for the complementary design.

\subsection{The Calderbank and Hobart inequalities}

\begin{Proposition}\label{Calderbank-ineq}
{\rm (\scaps{Calderbank} \cite{Calderbank88a})}
Let $\bar{x} = k-x$ and $\bar{y} = k-y$. Then
$$
(v-1)(v-2)\bar{x}\bar{y} - k(v-k)(v-2)(\bar{x}+\bar{y}) + k(v-k)(k(v-k)-1) \ge 0,
\eqno{(C)}
$$
with equality if and only if the design is a $3$-design. \qed
\end{Proposition}

Clearly, inequality (C) for a design is equivalent to this
inequality for the complementary design.
Calderbank observes that (C) is equivalent to (N).

\medskip
The following inequality was derived by Hobart
as a consequence of inequalities for coherent configurations.

\begin{Proposition}\label{SAH}
{\rm (\scaps{Hobart} \cite{Hobart95})}
The parameters of a quasisymmetric $2$-$(v,k,\lambda)$ design
with intersection numbers $x,y$, where $k > x > y$,
with strongly regular intersection-$x$ graph with eigenvalues
$K,R,S$, where $K > R > S$, satisfy
$$
\frac{v-2}{v} \left( 1 + \frac{R^3}{K^2} - \frac{(R+1)^3}{(b-K-1)^2} \right)
 - \frac{(v-2k)^2 \lambda}{k^2 (k-1)(v-k)} \ge 0.
\eqno{(H)}
$$
\end{Proposition}

This can also be formulated as
$Q^1_{11} \ge \frac{(v-2k)^2(v-1)}{k(v-k)(v-2)}$, where
$Q^1_{11}$ is the obvious Krein parameter of the strongly
regular graph.

Since the strongly regular graph (for the largest intersection size)
is the same for a quasisymmetric design and the complementary design,
we see that inequality (H) for a design is equivalent to this
inequality for the complementary design.

In the next section we show the equivalence of (C) and (H).

\section{Proof of Hobart's inequality}
Let $A = 1 + \frac{R^3}{K^2} - \frac{(R+1)^3}{(b-K-1)^2}$
be the parenthetical part of the inequality (H).
Substitute $b = V$ and $V = \frac{(K-R)(K-S)}{M}$ and $M = K+RS$ to get\\
$A = - \frac{(K-R)(KR+R^2-2KS+2R^2S-KS^2-RS^2)}{K^2(S+1)^2}$.
Now (H) says
{\small
$$
-\frac{v-2}{v} \, \frac{(K-R)(KR+R^2-2KS+2R^2S-KS^2-RS^2)}{K^2(S+1)^2}
- \frac{(v-2k)^2 \lambda}{k^2 (k-1)(v-k)} \ge 0.
$$}

If $S = -1$, then $x = k$ and the design is a multiple of a square
(or symmetric) design, a case that was excluded.
Hence $S < -1$. Multiply by $v K^2 (S+1)^2$ and substitute
$R = \frac{r-\lambda-k+y}{x-y}$ and $S = -\frac{k-y}{x-y}$
and $K = \frac{(r-1)k-(b-1)y}{x-y}$ and multiply by $(x-y)^4$
and substitute $\lambda = \frac{r(k-1)}{v-1}$ and $r = \frac{bk}{v}$
and multiply by $\frac{(v-1)^3}{b^3}$ and substitute
the value of $b$ found from equality in Proposition \ref{CCowen-ineq}.
Since we have $e > 0$ in Proposition \ref{nexus}, it follows that
$k \lambda \ne ry$, that is, $k^2-k-vy+y \ne 0$.
Divide by $(k^2-k-vy+y)^2$.
We see that (H) says
{\footnotesize
$$
(v-1)(v-2)xy+k^2(k-1)(k-3)+2k(k-1)(x+y) - k(k-1)v(x+y-1)
\ge 0
$$\par}

\noindent
but this is precisely inequality (C).

\medskip
In the same way one sees that Calderbank's inequality (C)
is equivalent to Neumaier's inequality (N).

\section{On the Blokhuis-Calderbank conditions}
Additional nonexistence results were given by
\scaps{Bagchi} \cite{Bagchi92} and
\scaps{Blokhuis} \& \scaps{Calderbank} \cite{BlokhuisCalderbank92}.
The methods and results are rather similar, but the results are
not equivalent: the latter paper eliminates several parameter sets
that survive other tests.
We do not repeat their definitions and results, but add some comments.
This is the table from \cite{BlokhuisCalderbank92}, p.~203.

\smallskip
\newcommand{\centercell}[1]{\multicolumn{1}{c}{#1}}
\begin{tabular}{rrrrrl}
\centercell{$v$} & \centercell{$k$} & \centercell{$\lambda$} &
\centercell{$y$} & \centercell{$x$} & comment \\
\hline
1090 & 540 & 2646 & 243 & 270 & fails \cite{BlokhuisCalderbank92}, Theorem~5.1 \\
1101 & 495 & 2223 & 198 & 225 &  \\
1266 & 396 & 1422 & 99 & 126 & fails \cite{BlokhuisCalderbank92}, Lemma 5.5 \\
1443 & 624 & 2136 & 246 & 273 & fails \cite{BlokhuisCalderbank92}, Theorem~5.1 \\
2704 & 544 & 1086 & 85 & 112 &  \\
2976 & 528 & 1023 & 69 & 96 & fails \cite{BlokhuisCalderbank92}, Theorem~5.1 for complement \\
5292 & 378 & 29 & 0 & 27 & fails \cite{Shrikhande53}, Theorem 3
\end{tabular}

\smallskip
In \cite{BlokhuisCalderbank92} it is said that Theorem 5.1 summarizes
the earlier results, but that theorem does not rule out the third
parameter set, while Lemma 5.5 does (but the theorem rules out
the complementary parameter set).

The last parameter set here is that of an $ARD(14,2)$,
where an affine resolvable design $ARD(n,t)$ is a 2-$(v,k,\lambda)$ design
with parameters $v = nk = n^2((n-1)t+1)$, $b = nr = n(n^2t+n+1)$,
$\lambda = nt+1$ where there is a resolution into $r$ parallel classes,
and any two blocks from different classes have $k^2/v = (n-1)t+1$
points in common. Using the Hasse invariant
\scaps{Shrikhande}~\cite{Shrikhande53} shows that no $ARD(n,t)$ exists
when $n \equiv 2$ (mod 4) and the square-free part of $n$ contains
a prime $\equiv 3$ (mod 4).

On the other hand, several far smaller parameter sets are ruled out.

\smallskip
\begin{tabular}{rrrrrrrl}
\centercell{$v$} & \centercell{$k$} & \centercell{$\lambda$} &
\centercell{$y$} & \centercell{$x$} & \centercell{$r$} & 
\centercell{$b$} & comment \\
\hline
77 & 33 & 24 & 12 & 15 & 57 & 133 & fails \cite{Bagchi92} and \cite{BlokhuisCalderbank92} \\
101 & 21 & 21 &  3 &  6 & 105 & 505 & fails \cite{Bagchi92} and \cite{BlokhuisCalderbank92} \\
137 & 40 &195 & 10 & 15 & 680 & 2329 & fails \cite{Bagchi92} and \cite{BlokhuisCalderbank92} \\
145 & 70 &161 & 28 & 35 & 336 & 696 & fails \cite{Bagchi92} and \cite{BlokhuisCalderbank92} \\
163 & 64 &672 & 22 & 28 & 1728 & 4401 & fails \cite{BlokhuisCalderbank92} \\
172 & 28 & 63 &  4 & 10 & 399 & 2451 & fails \cite{BlokhuisCalderbank92} \\
176 & 50 & 49 &  8 & 15  & 175 & 616 & fails \cite{BlokhuisCalderbank92} \\
\end{tabular}

\smallskip
In the first four cases, the complementary design violates \cite{Bagchi92}, Theorem~1.

\end{document}